\RequirePackage{ifpdf}
\ifpdf 
\documentclass[pdftex]{sigma}
\else
\documentclass{sigma}
\fi

\numberwithin{equation}{section}
\newtheorem{Theorem}{Theorem}[section]
\newtheorem{Proposition}[Theorem]{Proposition}
\newtheorem{Lemma}[Theorem]{Lemma}
\newtheorem{Corollary}[Theorem]{Corollary}
{\theoremstyle{definition}
\newtheorem{Definition}[Theorem]{Definition}
\newtheorem{Remark}[Theorem]{Remark}
}

\usepackage{mathrsfs}

\begin{document}


\renewcommand{\PaperNumber}{020}

\FirstPageHeading

\ShortArticleName{Geodesic Reduction via Frame Bundle Geometry}

\ArticleName{Geodesic Reduction via Frame Bundle Geometry}

\Author{Ajit BHAND}

\AuthorNameForHeading{A. Bhand}

\Address{Department of Mathematics, University of Oklahoma, Norman, OK, USA}
\Email{\href{mailto:abhand@math.ou.edu}{abhand@math.ou.edu}}

\ArticleDates{Received October 12, 2009, in f\/inal form February 18, 2010;  Published online February 22, 2010}

\Abstract{A manifold with an arbitrary af\/f\/ine connection is considered and the geodesic spray associated with the connection is studied in the presence of a Lie group action. In particular, results are obtained that provide insight into the structure of the reduced dynamics associated with the given invariant af\/f\/ine connection. The geometry of the frame bundle of the given manifold is used to provide an intrinsic description of the geodesic spray. A fundamental relationship between the geodesic spray, the tangent lift and the vertical lift of the symmetric product is obtained, which provides a key to understanding reduction in this formulation.}


\Classification{53B05; 53C05; 53C22; 58D19}

\section{Introduction}\label{section1}

The geometry of systems with symmetry has been an active area of research in the last several years. The study of manifolds with certain special geometric structure invariant under a Lie group action leads to what is known as reduction theory. Such questions arise in, for example, geometric mechanics. In this framework, the presence of symmetry allows the dynamics on a manifold to be studied on a lower dimensional manifold. In mechanics, there are at least three dif\/ferent ways of describing dynamics on a manifold, corresponding to the Lagrangian, Hamiltonian and af\/f\/ine connection formulations respectively. While the reduction theory for Lagrangian and Hamiltonian systems has been well developed (see \cite{AbMa:78, CeMaRa:01b, CeMaRa:01a, MaWe:74}), these results have been obtained by using variational analysis and symplectic geometry respectively. The main reason behind following this approach is the fact that the dynamics for such systems arises from variational principles which are manifested by symplectic structures in the Hamiltonian framework. However, when the dynamics on a manifold are given in terms of the geodesic equation of an af\/f\/ine connection, we cannot use variational analysis unless additional structure, such as a metric, is provided. Mechanical systems for which the dynamics are given by the geodesics of an af\/f\/ine connection that is not Levi-Civita include systems subjected to velocity constraints (see, for example, \cite{ADL:98, LM:97} and Section~\ref{section5.4}).

We consider an arbitrary af\/f\/ine connection on a manifold invariant under the action of a Lie group and provide results that enable us to decompose the reduced geodesic spray corresponding to the af\/f\/ine connection using tools from dif\/ferential geometry only. In other words, we do not use variational methods.
In arriving at our results, we come to a deeper understanding of the geometry of bundle of linear frames and its relationship with the geometry of the tangent bundle of the given manifold.

\looseness=1
The setup we consider is the following. Let $M$ be a manifold and $G$ a Lie group which acts on~$M$ in such a manner that $M$ is the total space of a principal bundle over~$M/G$. The Lie group $G$ also acts on the bundle $L(M)$ of linear frames over~$M$ via the lifted action.
It is known that there is a one-to-one correspondence between principal connections on $L(M)$ and af\/f\/ine connections on~$M$~\cite{kn:78}.
Let $\omega$ be a $G$-invariant principal connection on $L(M)$ with~$\nabla$ the corresponding af\/f\/ine connection on $M$. The geodesic spray $Z$ corresponding to $\nabla$ is a~\emph{second-order vector field} on the tangent bundle $TM$ with the pro\-perty that the projection of its integral curves correspond to geodesics on~$M$. Thus, to understand how the geodesics evolves under symmetry, $Z$ is the appropriate object to study. Since additional structure is not available, we exploit the geometry of the linear frame bundle in order to fully understand the meaning of the geodesic spray (which is classically def\/ined in local coordinates). The f\/irst signif\/icant step in this direction is to provide an intrinsic def\/inition of the geodesic spray that uses frame bundle geometry. We are able to provide such a def\/inition.

Moving ahead, we give a new interpretation of the geodesic invariance of a distribution on the manifold $M$ using frame bundle geometry and provide a new proof of a characterization due to Lewis \cite{ADL:98} using the symmetric product.

\looseness=1
Next, we turn our attention to understanding the reduced geodesic spray of a given con\-nection. Our main idea is that it is possible to study reduction using only geometric data and without variational analysis. To our knowledge, the proposed approach of using frame bundle geometry to study reduction is new. Even though the reduction method is discussed in the context of the geodesic spray, it can be applied to a general invariant second-order vector f\/ield.

The reduction procedure presented in this paper is based on the author's thesis \cite{Bhand:07}. Bullo and Lewis did some preliminary work in this direction~\cite{BL:07} and recently, Crampin and Mestdag~\cite{crme:09} have presented an approach which is similar in spirit to ours. They consider reduction and reconstruction of general second-order systems and provide a decomposition of the reduced system into three parts.

The def\/inition of the geodesic spray using frame bundle geometry, provided in Section~\ref{section3}, enables us to f\/ind a formula relating the geodesic spray and the vertical and complete lifts. In Section~\ref{section4} we present our reduction methodology along with the main result (Theorem~\ref{main}) which provides a new coordinate-free way of decomposing the reduced geodesic spray in terms of objects def\/ined on the reduced space. Finally, in Section~\ref{section5} we provide a geometric interpretation of the result in the Riemannian case and discuss some avenues for future work.

\section{Linear connections}\label{section2}

In this section we review some concepts relevant to our investigation and establish notation to be used throughout the paper.
Let $M$ be an $n$-dimensional smooth manifold and $L(M)(M{,} GL(n{;}\mathbb{R}))\!$ the bundle of linear frames with total space $L(M)$, base space $M$, structure group $GL(n;\mathbb{R})$ and canonical projection $\pi_M$. We denote the (right) action of $GL(n;\mathbb{R})$ by $\Phi: L(M) \times GL(n;\mathbb{R}) \rightarrow L(M)$.  For f\/ixed $a \in GL(n;\mathbb{R})$, this action induces a map $\Phi_a: L(M) \rightarrow L(M)$ given by $\Phi_a(u)= \Phi(u, a)$.
Recall that a linear frame $u$ at $x\in M$ is an ordered basis $(X_1, X_2, \ldots, X_n)$ of the tangent space $T_xM$. If $a= (a^i_j) \in GL(n;\mathbb{R})$ and $u= (X_1,X_2, \ldots, X_n)$, then $ua:= \Phi_a(u)$ is the linear frame $(Y_1, Y_2, \ldots, Y_n)$ at $x$ def\/ined by $Y_i= \sum^n_{j=1}a^j_iX_j$. Equivalently, if $(e_1, e_2, \ldots, e_n)$ is a standard basis for $\mathbb{R}^n$,  a frame $u= (X_1, X_2, \ldots, X_n)$ at $x$ can also be def\/ined as a linear isomorphism $u:\mathbb{R}^n\rightarrow T_xM$ given by $ue_i= X_i$. In other words, if $\xi= \sum^n_{i=1} \xi^ie_i \in \mathbb{R}^n$, then $u\xi= \sum^n_{i=1}\xi^iX_i$. This is the notion of a linear frame that we will use throughout the paper.
 The inf\/initesimal generator corresponding to an element  $A \in \mathfrak{gl}(n,\mathbb{R})$ of the Lie algebra will be denoted by $A_{L(M)}$.
The canonical form $\theta$ of $L(M)$ is a one-form on $L(M)$ def\/ined by
$\theta (X_u)= u^{-1}\left(\pi_M(X_u)\right),\; X_u \in T_uL(M)$.

A principal connection $\omega$ in the bundle $L(M)(M, GL(n;\mathbb{R}))$ of linear
frames over $M$ is called a~\emph{linear connection on $M$}. The associated horizontal subbundle is denoted by $HL(M)$. The tangent bundle of $M$, denoted by $\tau_M{:}\,TM\! \rightarrow\! M$, is a bundle associated with $L(M)(M{,} GL(n{;}\mathbb{R}))$~\cite{kn:78}.

Given a linear connection $\omega$ on $M$, for each $\xi \in \mathbb{R}^n$, the \emph{standard horizontal vector field corresponding to $\xi$}, denoted by $B(\xi):L(M)\rightarrow TL(M)$, is def\/ined as follows. For each $u \in L(M)$, the vector $B(\xi)_u$ is the unique horizontal vector at $u$ with the
property that $T_u\pi_M (B(\xi)_u)= u\xi$.

We consider an arbitrary principal bundle $P(M,G)$ with total space $P$, base space $M$ and structure group $G$. The canonical vertical bundle will be denoted by $VP$.

Given a principal f\/iber bundle $P(M, G)$ and a representation
$\rho$ of $G$ on a f\/inite-dimensional vector space $V$, a
\emph{pseudotensorial $r$-form of type $(\rho, V)$ on $P$} is a
$V$-valued $r$-form $\varphi$ on $P$ such that
\begin{gather*}
\Phi^{\ast}_g\varphi= \rho\big(g^{-1}\big)\cdot \varphi, \qquad g \in G,
\end{gather*}
where $\Phi$ is the action of $G$ on $P$. A pseudotensorial $r$-form $\varphi$ of type $(\rho, V)$ is called
a \emph{tensorial $r$-form} if it is horizontal in the sense that
$\varphi(X_1,
\ldots,
X_r)=0$ whenever $X_i$ is vertical for at least one $i \in \{1,\ldots, r\}$.
A connection one-form $\omega$ on a principal bundle $P(M, G)$ is a pseudotensorial one-form of type $(\mathrm{Ad}(G), \mathfrak{g})$, where $\mathrm{Ad}(G)$ is the adjoint representation of $G$ on $\mathfrak{g}$.
The following result characterizing the set of all principal connections on $P$ can be readily proved.

\begin{Proposition}
\label{setofpconn}
Let $\omega$ be a principal connection one-form on a principal bundle $P(M, G)$ and let $\alpha$ be a tensorial one-form of type $(\mathrm{Ad}(G), \mathfrak{g})$ on $P$. Then $\bar{\omega}:= \omega+ \alpha$ defines a new principal connection on $P$. Conversely, given any two principal connection forms $\omega$ and $\bar{\omega}$ respectively, the object $\alpha:= \bar{\omega}- \omega$ is a tensorial one-form of type $(\mathrm{Ad}(G), \mathfrak{g})$ on $P$.
\end{Proposition}

In other words, the space of principal connections is an af\/f\/ine space modeled on the vector space of tensorial one-forms of type $(\mathrm{Ad}(G),\mathfrak{g})$.

To each vector f\/ield $Y$ on $M$ we can associate a function $f_Y: L(M)\rightarrow \mathbb{R}$ as follows. For $u \in L(M)$, we have $f_Y(u)= u ^{-1}(Y(\pi_M(u)))$.
The following result provides a correspondence between tensorial one-forms on $L(M)$ and $(1, 2)$ tensor f\/ields on $M$.
\begin{Proposition}
There is a one-to-one correspondence between tensorial one-forms of type $(\mathrm{Ad}(GL(n; \mathbb{R})), \mathfrak{gl}(n; \mathbb{R}))$ on $L(M)$ and $(1, 2)$ tensor fields on $M$.
\end{Proposition}
\begin{proof}
Since the tangent bundle $\tau_M:TM\rightarrow M$ is the bundle associated with $L(M)$ with standard f\/iber $\mathbb{R}^n$, for each $(1,2)$ tensor f\/ield $S$ on $M$, and $u \in L(M)$, we can def\/ine a map $\alpha_S: TL(M) \rightarrow \mathfrak{gl}(n; \mathbb{R})$ as follows. Let $\widetilde{X} \in T_uL(M)$, for $u \in L(M)$ and let $\eta \in \mathbb{R}^n$. Then
\begin{gather}
\label{eqref:alphas}
u\big(\alpha_S(u)(\widetilde{X})\eta\big)= S\big(T_u\pi_M(\widetilde{X}), u\eta\big).
\end{gather}
Since $\alpha_S(\widetilde{X}) \in \mathfrak{gl}(n; \mathbb{R})$, the product $\alpha(\widetilde{X})\eta \in \mathbb{R}^n$. We now show that $\alpha_S$ is a tensorial form of type $(\mathrm{Ad}(GL(n; \mathbb{R})), \mathfrak{gl}(n; \mathbb{R}))$. For $a \in GL(n; \mathbb{R})$, let $\widetilde{Y}= T_u\Phi_a\widetilde{X} \in T_{ua}L(M)$. Then, using the def\/inition \eqref{eqref:alphas}, we get
\begin{gather*}
(ua)\big(\alpha_S(ua)(\widetilde{Y})\eta\big)= S\big(T_{ua}\pi_M(\widetilde{Y}), ua\eta\big)= u\big(\alpha_S(u)(\widetilde{X})(a\eta)\big)
\end{gather*}
from which we get
\begin{gather*}
\alpha_S(ua)(\widetilde{Y})= a^{-1}\alpha_S(u)(\widetilde{X})a,
\end{gather*}
which means that $\alpha_S$ is pseudotensorial. Next, if $\widetilde{X} \in V_uL(M)$, it is easy to see that $\alpha_S(\widetilde{X})= 0$, which shows that $\alpha_S$ is tensorial.

Conversely, given a tensorial one-form $\alpha:TP \rightarrow \mathfrak{gl}(n; \mathbb{R})$, we can def\/ine a $(1, 2)$ tensor f\/ield~$S_{\alpha}$ as follows:
\begin{gather*}
S_{\alpha}(X, Y)= u\big(\alpha(\widetilde{X})f_Y(u)\big),\qquad X, Y \in T_xM,\qquad  \pi_M(u)= x,
\end{gather*}
where $\widetilde{X}_u \in T_uL(M)$ has the property that $T_u\pi_M(\widetilde{X})= X$. Since $\alpha$ is tensorial, $S_{\alpha}$ is well-def\/ined.
\end{proof}
\begin{Corollary}
Let $\omega$ and $\bar{\omega}$ be linear connections of $M$ and let $\nabla$ and $\overline{\nabla}$, respectively, be the corresponding covariant derivatives. If $\alpha= \bar{\omega}- \omega$ then, for vector fields $X$ and $Y$ on $M$, we have
\begin{gather*}
\overline{\nabla}_XY= \nabla_XY+ S_{\alpha}(X, Y),
\end{gather*}
where $S_{\alpha}$ is the unique $(1,2)$ tensor field on $M$ corresponding to $\alpha$.
\end{Corollary}

This result, therefore, characterizes the set of all af\/f\/ine connections on the manifold $M$.

\section{The geodesic spray of an af\/f\/ine connection}\label{section3}

In this section we study the geodesic spray associated with a given af\/f\/ine connection. This object is typically def\/ined in terms of local coordinates and here we provide an intrinsic def\/inition using the geometry of the linear frame bundle.
Given a linear connection $\omega$ on $M$, for \emph{fixed} $\xi \in \mathbb{R}^n$, let $\Phi_{\xi}: L(M) \rightarrow TM$ be the \emph{association map} given by
$\Phi_{\xi}(u)= u\xi$.
We def\/ine a (second-order) vector f\/ield $Z: TM \rightarrow TTM$ called the \emph{geodesic spray} as follows:
\begin{gather}
\label{eqref:geospray1}
Z(v)= T_u\Phi_{\xi}(B(\xi)_u), \qquad v \in TM,
\end{gather}
where $u \in L_{\tau_M(v)}(M)$ and $\xi \in \mathbb{R}^n$ are such that $u\xi= v$, and
$B(\xi)$ is the standard horizontal vector f\/ield corresponding to $\xi$ for the linear connection $\omega$ associated with $\nabla$.
We have the following result.
\begin{Proposition}
The map  $Z$ defined in \eqref{eqref:geospray1} is a second-order vector field on $TM$. The coordinate expression for $Z$, in terms of the canonical tangent bundle coordinates $(x^i, v^i)$ is given by
\begin{gather}
\label{eqref:geosprcoord1}
Z= v^i\frac{\partial}{\partial x^i}- \Gamma^i_{jk}v^jv^k\frac{\partial}{\partial v^i}.
\end{gather}
\end{Proposition}
\begin{proof}
We f\/irst show that $Z$ as given by \eqref{eqref:geospray1} is well-def\/ined. The canonical projection on the tangent bundle is denoted by $\tau_M: TM \rightarrow M$. For a given $v \in TM$, we write $x:= \tau_M(v)$. Suppose that $u^{\prime} \in L_x(M)$ \ and $\xi^{\prime} \in \mathbb{R}^n$ are such that $u^{\prime}\xi^{\prime}= v= u \xi$. Then, we must have $u^{\prime}= ua$ for some $a \in GL(n;\mathbb{R})$. Consequently, $\xi^{\prime}= a^{-1}\xi$. We compute
\begin{gather*}
T_{ua}\Phi_{a^{-1}\xi}\big(B(a^{-1}\xi)_{ua}\big) =  T_{ua}\Phi_{a^{-1}\xi} T_uR_a(B(\xi)_u)
 = T_u(\Phi_{a^{-1}\xi}\circ R_a)B(\xi)_u = T_u\Phi_{\xi}(B(\xi)_u),
\end{gather*}
where the f\/irst equality follows from the properties of a standard horizontal vector f\/ield. Let us now show that $Z$ is a second-order vector f\/ield. We have
\begin{gather*}
T\tau_M(Z(v)) = T\tau_M(T_u\Phi_{\xi}B(\xi)_u)= T_u(\tau_M \circ \Phi_{\xi})B(\xi)_u
 = T_u(\pi_M)(B(\xi)_u)= u\xi=  v
\end{gather*}
as desired.

It now remains to be shown that the coordinate representation of $Z$ is as given in \eqref{eqref:geosprcoord1}, but this follow directly from the coordinate representation of $B(\xi)$.
\end{proof}

\subsection{Tangent and vertical lifts}\label{section3.1}

If $X$ is a vector f\/ield on $M$ we can def\/ine a unique vector f\/ield $\widetilde{X}$ on $L(M)$ corresponding to $X$ as follows. Let $\phi_t^X$ be the f\/low of $X$. The \emph{tangent lift} $X^T$ is a vector f\/ield on $TM$ def\/ined by
\begin{gather*}
X^T(v_x)= \left.\frac{d}{dt}\right|_{t= 0} T\phi_t^X(v_x).
\end{gather*}
Let $u \in L_x(M)$ and $\xi \in \mathbb{R}^n$ be such that $u\xi= v_x$. For $\xi$ f\/ixed, let $\Phi_{\xi}: L(M) \rightarrow TM$ be the association map. The f\/low of $X^T$ def\/ines a curve $u_t$ in $L(M)$ by $u_t=  T_x\phi_t^X \cdot u$. That is,
\begin{gather*}
\Phi_{\xi}u_t= \big(T_x\phi_t^X\circ \Phi_{\xi}\big)u.
\end{gather*}
The map $\widetilde{\Phi}_t(u)= u_t$ def\/ines a f\/low on $L(M)$. The corresponding vector f\/ield is called the \emph{natural lift} $\widetilde{X}$ of $X$ onto $L(M)$. Thus, we have
\begin{gather*}
X^T(v_x)= T_u\Phi_{\xi}\widetilde{X}(u).
\end{gather*}
Given $v_x, w_x \in T_xM$, the \emph{vertical lift} of $w$ at $v$ is def\/ined by
\begin{gather*}
\mathrm{vlft}_{v_x}(w_x)= \left.\frac{d}{dt}\right|_{t= 0} (v_x+ tw_x).
\end{gather*}

\subsection{Decomposition of the geodesic spray}\label{section3.2}

The following result provides an explicit relationship between $\nabla$ and the connection one-form $\omega$ of the corresponding linear connection.
\begin{Proposition}
\label{sympomega}
Let $M$ be a manifold with a connection $\nabla$ with the corresponding linear connection one-form $\omega$. Given vector fields $X$ and $Y$ on $M$, let $\widetilde{X}$ be the natural lift of $X$ onto $L(M)$, and $f_{Y}: L(M) \rightarrow \mathbb{R}^n$ the function associated with $Y$. Then,
\begin{gather}
\label{eqref:nablaomega}
\nabla_XY(x)= [X, Y](x)+ u\big(\omega(\widetilde{X}(u))f_Y(u)\big),\qquad \pi_M(u)= x.
\end{gather}
\end{Proposition}

\begin{proof}
Let us f\/irst verify that the right-hand side of \eqref{eqref:nablaomega} is independent of the choice of $u \in L_x(M)$. For $a \in GL(n; \mathbb{R})$, we compute
\begin{gather*}
(ua)\big(\omega(\widetilde{X}_{ua})f_Y(ua)\big)  = (ua)\big(\omega(T_u\Phi_a\widetilde{X}_u)a^{-1}f_Y(u)\big)\\
\phantom{(ua)\big(\omega(\widetilde{X}_{ua})f_Y(ua)\big)}{}  = (ua)\big(a^{-1}\omega(\widetilde{X}_u)a (a^{-1}f_Y(u))\big)
 = u\big(\omega(\widetilde{X}_u)f_Y(u)\big).
\end{gather*}

We shall now prove that
\begin{gather*}
\mathscr{L}_{X^h}f_Y(u)= u^{-1}\left([X, Y](x)\right)+ \omega\big(\widetilde{X}(u)\big)f_Y(u).
\end{gather*}
Notice that $\mathrm{hor}(\widetilde{X})= X^h$ since both are horizontal vector f\/ields on $L(M)$ projecting to~$X$. In other words,
\begin{gather*}
X^h= \widetilde{X}- (\omega(\widetilde{X}))_{L(M)},
\end{gather*}
where $\omega(\widetilde{X})_{L(M)}$ is the vertical vector f\/ield on $L(M)$ given by
\begin{gather*}
u \mapsto \omega(\widetilde{X}(u))_{L(M)}(u).
\end{gather*}
Therefore,
\begin{gather}
\label{eq:lied2}
\mathscr{L}_{X^h}f_{Y}(u)= \mathscr{L}_{\widetilde{X}}f_Y(u)- \mathscr{L}_{\omega(\widetilde{X})_{L(M)}}f_Y(u).
\end{gather}
The f\/low of $\widetilde{X}$ is $\Phi^{\widetilde{X}}_t(u)= T_x\Phi^X_t\cdot u$, where $\Phi^X_t$ is the f\/low of $X$ and $x=\pi_M(u)$.
The f\/irst term on the right-hand side of \eqref{eq:lied2} is
\begin{gather*}
\mathscr{L}_{\widetilde{X}}f_Y(u)  = \left. \frac{d}{dt}f_Y\big(T_x\Phi^X_t\cdot u\big)\right|_{t= 0}
 = \left. \frac{d}{dt}\big(u^{-1}T_{\Phi^X_t(x)}\Phi^X_{-t}\big)\big(Y(\Phi^X_t(x))\big)\right|_{t= 0}\\
\phantom{\mathscr{L}_{\widetilde{X}}f_Y(u)}{}
= \left.\frac{d}{dt}u^{-1}{\Phi^x_t}^{\ast}Y(x)\right|_{t= 0}= u^{-1}[X, Y](x).
\end{gather*}
We next compute the second term on the right-hand side of \eqref{eq:lied2}
\begin{gather*}
\mathscr{L}_{\omega(\widetilde{X})_{L(M)}}f_Y(u)
 = \left.\frac{d}{dt}f_Y\big(u\exp \big(t \omega(\widetilde{X}(u))\big)\big)\right|_{t= 0}
= \left. \frac{d}{dt}\big(\exp \big(-t\omega(\widetilde{X}(u))\big)u^{-1}Y(x)\big)\right|_{t= 0} \\
\phantom{\mathscr{L}_{\omega(\widetilde{X})_{L(M)}}f_Y(u)}{}
 = -\omega(\widetilde{X}(u))u^{-1}Y(x)= -\omega(\widetilde{X}(u))f_Y(u),
\end{gather*}
where the second equality above follows from the fact that $f_X$ is a pseudotensorial form of degree zero. This completes the proof.
\end{proof}

\begin{Remark}
It is known that every derivation $\mathscr{D}$ of the tensor algebra of $M$ can be decomposed~as
\begin{gather*}
\mathscr{D}= \mathscr{L}_X+ S,
\end{gather*}
where $X$ is a vector f\/ield on $M$ and $S$ is a $(1, 1)$ tensor f\/ield on $M$.
From the previous result, it follows that the $(1, 1)$ tensor f\/ield associated with the derivation $\nabla_X$ is given by
\begin{gather*}
T_xM \ni v \mapsto u\big(\omega(\widetilde{X}(u))f_{Y_v}(u)\big),\qquad \pi_M(u)= x,
\end{gather*}
where $Y_v$ is any vector f\/ield on $M$ with value $v$ at $x$.
\end{Remark}

We now prove the main result of this section.
\begin{Proposition}\label{zdecomp}
Let $v \in T_xM$ for some $x \in M$, and $X_v$ be an arbitrary vector field that has the value $v$ at $x$. Then,
\begin{gather}
\label{eqref:zxtvlf}
Z(v)= (X_v)^T(v)- \mathrm{vlft}_v(\nabla_{X_v}X_v(x)).
\end{gather}
\end{Proposition}
\begin{proof}
Using local coordinates $x^i$ around $x$ in $M$, we write $X_v= {X_v}^i\frac{\partial}{\partial x^i}$. Then,
\begin{gather*}
(X_v)^T(v)= v^i\frac{\partial}{\partial x^i}+ v^j\frac{\partial {X_v}^i}{\partial x^j}\frac{\partial}{\partial v^i},
\end{gather*}
and
\begin{gather*}
\nabla_{X_v}X_v(x)= \left(\frac{\partial {X_v}^i}{\partial x^j}{X_v}^j+
\Gamma^i_{jk}{X_v}^j{X_v}^k\right) \frac{\partial}{\partial x^i}=  \left(\frac{\partial {X_v}^i}{\partial x^j}v^j+ \Gamma^i_{jk}v^jv^k\right) \frac{\partial}{\partial x^i} .
\end{gather*}
So
\begin{gather*}
\mathrm{vlft}_v(\nabla_{X_v}X_v(x))=  \left(\frac{\partial {X_v}^i}{\partial x^j}v^j+ \Gamma^i_{jk}v^jv^k\right) \frac{\partial}{\partial v^i}.
\end{gather*}
Thus,
\begin{gather*}
(X_v)^T(v)- \mathrm{vlft}_v(\nabla_{X_v}X_v(x)) =   v^i\frac{\partial}{\partial x^i}+ v^j\frac{\partial {X_v}^i}{\partial x^j}\frac{\partial}{\partial v^i}-
\left(\frac{\partial {X_v}^i}{\partial x^j}v^j+ \Gamma^i_{jk}v^jv^k\right) \frac{\partial}{\partial v^i}\\
 \phantom{(X_v)^T(v)- \mathrm{vlft}_v(\nabla_{X_v}X_v(x))}{}
 = v^i\frac{\partial}{\partial x^i}- \Gamma^i_{jk}v^jv^k \frac{\partial}{\partial v^i}= Z(v).
\end{gather*}
This proves the result.  Notice that, even though each of the two terms $X_v^T(v)$ and $\mathrm{vlft}_v(\nabla_{X_v}X_v)$ depends on the extension $X_v$, the terms that depend on the derivative of $X_v$ cancel in the expression for $Z$.

\noindent {\bf Alternate proof.} Given $v_x \in T_xM$, let $u \in L_xM$ and $\xi \in \mathbb{R}^n$ be such that $u\xi= v_x$. Let $X_v$ be a vector f\/ield with value $v_x$ at $x$ and $f_{X_v}: L(M) \rightarrow \mathbb{R}^n$ be the corresponding function on $L(M)$.
Then, for each $u \in L(M)$,
\begin{gather*}
\mathrm{hor}(\widetilde{X_v}(u))= B(f_{X_v}(u))_u,
\end{gather*}
where $\widetilde{X_v}$ is the natural lift of $X_v$ onto $L(M)$. We have
\begin{gather*}
Z(X_v(x))  = T_u\Phi_{f_{X_v}(u)}B(f_{X_v}(u))_u
  =  T_u\Phi_{f_{X_v}(u)}\big(\widetilde{X}(u)- \omega(\widetilde{X}(u))_{L(M)}\big)\\
\phantom{Z(X_v(x))}{}  = X^T_v(X_v(x))-  T_u\Phi_{f_{X_v}(u)}\big(\omega(\widetilde{X}(u))_{L(M)}\big)\\
\phantom{Z(X_v(x))}{} = X^T_v(X_v(x))- \left. \frac{d}{dt}\Phi_{f_{X_v}(u)}u\exp\big(t\omega(\widetilde{X}(u))\big)\right|_{t= 0}\\
\phantom{Z(X_v(x))}{} =  X^T_v(X_v(x))- \mathrm{vlft}_{X_v(x)}u\big(\omega(\widetilde{X}(u))f_{X_v}(u)\big)\\
\phantom{Z(X_v(x))}{} =  X^T_v(X_v(x))- \mathrm{vlft}_{X_v(x)}\nabla_{X_v}X_v(x),
\end{gather*}
where we have used Proposition~\ref{sympomega} in the last step.
\end{proof}

\subsection{Geodesic invariance}\label{section3.3}

We recall the notion of geodesic invariance.
\begin{Definition}
A distribution $D$ on a manifold $M$ with an af\/f\/ine connection $\nabla$ is called \emph{geodesically invariant} if for every geodesic $c:[a, b]\rightarrow M$,
$\dot{c}(a) \in D_{c(a)}$ implies that $\dot{c}(t) \in D_{c(t)}$ for all $t \in [a, b]$.
\end{Definition}

It turns out that geodesic invariance can be characterized by studying a certain product on the set of vector f\/ields on $M$. Let $M$ be a manifold with a connection $\nabla$. Given vector f\/ields $X, Y \in \Gamma(TM)$, the \emph{symmetric product} $\langle X: Y\rangle $ is the vector f\/ield def\/ined by
\begin{gather}
\label{eqref:sympro}
\langle X: Y \rangle = \nabla_XY+ \nabla_YX.
\end{gather}

The following result gives a description of the symmetric product using linear frame bundle geometry.
\begin{Proposition}
\label{prop:sympframe1}
Let $M$ be a manifold with a connection $\nabla$, let $\omega$ be the associated linear connection and $\theta$ the canonical form on $L(M)$ respectively. If $X$ and $Y$ are vector fields on $M$ and $\widetilde{X}$ and $\widetilde{Y}$ the respective natural lifts on $L(M)$, then
\begin{gather}
\label{sympframe}
\langle X: Y\rangle (x)= 2u\big(\mathrm{Sym}(\omega \otimes \theta)(\widetilde{X}, \widetilde{Y})\big), \qquad \pi_M(u)= x.\end{gather}
\end{Proposition}
\begin{proof}
It is clear that the right-hand side of \eqref{sympframe} is independent of the choice of $u \in L_x(M)$.
By def\/inition, we have $\theta(\widetilde{X}(u))= f_X(u)$ and $\theta(\widetilde{Y}(u))= f_Y(u)$.
We compute
\begin{gather*}
2u\big(\mathrm{Sym}(\omega\otimes \theta)(\widetilde{X}(u), \widetilde{Y}(u))\big)  = u\big(\omega(\widetilde{X}(u))\theta(\widetilde{Y}(u))+ \omega(\widetilde{Y}(u))\theta(\widetilde{X}(u))\big)\\
\phantom{2u\big(\mathrm{Sym}(\omega\otimes \theta)(\widetilde{X}(u), \widetilde{Y}(u))\big)}{}
= u\big(\omega(\widetilde{X}(u))f_Y(u)+ \omega(\widetilde{Y}(u))f_X(u)\big).
\end{gather*}
The result now follows from Proposition \ref{sympomega}.
\end{proof}

\begin{Remark} The object $\mathrm{sym}(\omega\otimes \theta)$ def\/ines a quadratic form $\Sigma_u:T_uL(M)\times T_uL(M) \rightarrow T_uL(M)$, $u \in L(M)$ as follows:
\begin{gather*}
\Sigma_u(X, Y)= \widetilde{Z_{X,Y}}(u),\qquad X, Y \in T_uL(M),
\end{gather*}
where $\widetilde{Z_{X,Y}}$ is the natural lift onto $L(M)$ of the vector f\/ield $Z_{X,Y}$ on $M$ given by
\begin{gather*}
Z_{X,Y}(x)= u(\mathrm{sym}(\omega \otimes \theta)(X, Y)),\qquad \pi(u)= x.
\end{gather*}
This is seen to be well-def\/ined.
\end{Remark}

Given a distribution $D$ on $M$, we represent by $\Gamma(D)$ the set of vector f\/ields taking values in~$D$.
The following result, proved by Lewis \cite{ADL:98}, provides inf\/initesimal tests for geodesic invariance and gives the geometric meaning of the symmetric product.
\begin{Theorem}[Lewis]
\label{thm:geoinv}
Let $D$ be a distribution on a manifold $M$ with a connection $\nabla$. The following are equivalent:
\begin{enumerate}\itemsep=0pt
\item[$(i)$] $D$ is geodesically invariant;
\item[$(ii)$] $\langle X: Y \rangle  \in \Gamma(D)$ for every $X, Y \in \Gamma(D)$;
\item[$(iii)$] $\nabla_XX \in \Gamma(D)$ for every $X \in \Gamma(D)$.
\end{enumerate}
\end{Theorem}

We give an intrinsic proof of this theorem below.
Thus, for geodesically invariant distributions, the symmetric product plays the role that the Lie bracket plays for integrable distributions. We use this result to interpret some terms obtained in the decomposition of the reduced geodesic spray in Section~\ref{section4}.

Now, given a $p$-dimensional distribution $D$ on an $n$-dimensional manifold $M$ with a linear connection, we say that a frame $u \in L_x(M)$ is \emph{$D$-adapted} if $u|_{\mathbb{R}^p}$ is an isomorphism onto $D_x$. Let $L(M, D)$ be the collection of $D$-adapted frames. We observe that $L(M, D)$ is invariant under the subgroup of $GL(n; \mathbb{R})$ consisting of those automorphisms which leave $\mathbb{R}^p$ invariant.  It turns out that $L(M, D)$ is a subbundle of $L(M)$  with structure group $H$ given by
\begin{gather*}
H= \left\{ A \in GL(n;\mathbb{R})\,|\, A= \left(\!\begin{array}{cc}a&b\\0&c\end{array}\!\right),\, a \in GL(p; \mathbb{R}),\, b \in L(\mathbb{R}^{n-p}, \mathbb{R}^p),\, c \in GL(n-p; \mathbb{R})\right\}.
\end{gather*}
We denote the bundle of $D$-adapted frames by $L(M, D)(M, H)$ and the Lie algebra of $H$ by $\mathfrak{h}$.
We have the following result.
\begin{Proposition}
The distribution $D$ is geodesically invariant if and only if, for each $\xi \in \mathbb{R}^p$, $B(\xi \oplus 0)|_{L(M, D)}$ is a vector field on $L(M, D)$.
\end{Proposition}
\begin{proof}
We f\/irst prove the ``if'' statement.  Suppose that $B(\xi\oplus 0)$ is a vector f\/ield on $L(M,D)$ and let $c: \mathbb{R} \rightarrow L(M)$ be its integral curve passing through $\bar{u} \in L(M, D)$. Then, we know that $x(t):= \pi_M(c(t))$ is the unique geodesic with the initial condition $\bar{u}\xi \in D$. We must show that $\dot{x}(t) \in D_{x(t)}$ for all $t$.
We have
\begin{gather*}
\dot{x}(t)= T\pi_M(B(\xi\oplus 0)_{c(t)})= c(t)(\xi\oplus 0).
\end{gather*}
Since $B(\xi \oplus 0)$ is a vector f\/ield on $L(M,D)$, we must have $c(t) \in L(M, D)$ for all $t$. Thus, we have $\dot{x}(t) \in D_{x(t)}$ for all $t$. The ``only if'' part of the statement can be proved by reversing this argument.
\end{proof}
An immediate consequence of this result is the following.
\begin{Corollary}
A distribution $D$ is geodesically invariant if and only if the geodesic spray $Z$ is tangent to the submanifold $D$ of $TM$.
\end{Corollary}
We are now in a position to provide a proof of Theorem \ref{thm:geoinv} using frame bundle geometry.
\begin{proof}[Proof of Theorem~\ref{thm:geoinv}]
$(i) \Longrightarrow  (ii)$ Suppose that $D$ is geodesically invariant, and let $X_1, X_2 \in \Gamma(D)$. Then, we know that the corresponding functions $f_{X_i}: L(M, D) \rightarrow \mathbb{R}^p \oplus \mathbb{R}^{n-p},\; i=1,2$, take values in $\mathbb{R}^p$. Also,
\begin{gather*}
(X_i)^h (u)= c^j_iB(e_j \oplus 0)_u, \qquad u \in L(M, D),
\end{gather*}
where $c^j_i$ are functions on $L(M)$ and $\{e_j\}_{j= 1,\ldots, p}$ is the standard basis for $\mathbb{R}^p$. This is possible since $\{B(e_i)\}$ form a basis for $H_uL(M)$.
We have
\begin{gather*}
f_{(\nabla_{X_1}X_2+ \nabla_{X_2}X_1)}= \mathscr{L}_{(X_1)^h}f_{X_2}+ \mathscr{L}_{(X_2)^h}f_{X_1} = c^j_1\mathscr{L}_{B(e_j\oplus 0)}f_{X_2}+ c^k_2 \mathscr{L}_{B(e_k\oplus 0)}f_{X_1}.
 \end{gather*}
 Since $f_{X_i}$, $i= 1, 2$, are $\mathbb{R}^p$-valued functions on $L(M, D)$ and $B(e_j\oplus 0)|_{L(M, D)}$, $j= 1,\ldots, p$, are vector f\/ields on $L(M, D)$ because the distribution is assumed to be geodesically invariant, we conclude that the function $f_{(\nabla_{X_1}X_2+ \nabla_{X_2}X_1)}: L(M,D) \rightarrow \mathbb{R}^p\oplus \mathbb{R}^{n-p}$ takes its values in $\mathbb{R}^p$. This proves $(ii)$.

$(ii) \Longrightarrow (iii)$ This follows directly from the def\/inition of the symmetric product.

$(iii) \Longrightarrow (i)$ Assume that $\nabla_XX \in \Gamma(D)$ for every $X \in \Gamma(D)$. This implies that the function $\mathscr{L}_{X^h}f_X: L(M, D) \rightarrow \mathbb{R}^p \oplus \mathbb{R}^{n- p}$ takes values in $\mathbb{R}^p$.   Once again, we can write $X^h= C^iB(e_i\oplus 0)$ for some functions $C^i$. This implies that $B(e_i\oplus 0)|_{L(M, D)}$ must be a vector f\/ield on $L(M, D)$.
\end{proof}

The above result shows that it is possible to check for geodesic invariance by looking at vector f\/ields $B(e_i\oplus 0)$ on the bundle $L(M,D)$.

Theorem \ref{thm:geoinv} and Proposition \ref{prop:sympframe1} suggest the following:
\begin{Corollary}
Let $D$ be a $p$-dimensional distribution on a manifold $M$ with a linear connection $\omega$. Let $\widetilde{D}$ be the natural lift of $D$ onto $L(M)$ and $L(M, D)(M, H)$ the bundle adapted to $D$. The following statements are equivalent:
\begin{enumerate}\itemsep=0pt
\item[$(i)$] $D$ is geodesically invariant;
\item[$(ii)$] $\mathrm{Sym}(\omega \otimes \theta)$ is an $\mathbb{R}^p$-valued quadratic form on $\widetilde{D}|_{L(M,D)}$;
\item[$(iii)$] $\omega(\widetilde{X}(u)) \in \mathfrak{h}\; \mathrm{for\; all\;} \widetilde{X} \in \Gamma(\widetilde{D})$, $u \in L(M,D)$.
\end{enumerate}
\end{Corollary}

\section{Geodesic reduction}\label{section4}

In this section we consider the following setup. Let $M(M/G, G)$ be a principal f\/iber bundle with structure group $G$ and let $\nabla$ be a $G$-invariant af\/f\/ine connection on $M$.
We study the ``reduced'' geodesic spray corresponding to $\nabla$  and use a principal connection to decompose it into various components. These components correspond to geometric objects def\/ined on $M/G$.

\subsection{Invariant af\/f\/ine connections}\label{section4.1}
Let $M(M/G, G)$ be a principal f\/iber bundle with a $G$-invariant af\/f\/ine connection $\nabla$ on $M$. Choose a principal connection $A$ on $M(M/G, G)$. With this data, we can def\/ine an af\/f\/ine connection~$\nabla^A$ on $M/G$ as follows.
\begin{Proposition}
Let $\nabla$ and $A$ be as above.
Given vector fields $X$ and $Y \in \Gamma(T(M/G))$ on~$M/G$, the map $\nabla^A: \Gamma(T(M/G))\times \Gamma(T(M/G))\rightarrow \Gamma(T(M/G))$ defined by
\begin{gather}
\label{eq:redaff}
\nabla^A_XY(x)= T\pi_{M/G}\nabla_{X^h}Y^h(q),\qquad q \in M,\qquad \pi_{M/G}(q)= x \in M/G
\end{gather}
is an affine connection on $M/G$.
\end{Proposition}
\begin{proof}
It is easy to see that $\nabla^A$ is well-def\/ined. Given a smooth function $f: M/G \rightarrow \mathbb{R}$, we def\/ine the lift $\tilde{f}:M\rightarrow \mathbb{R}$ by $\tilde{f}= \pi^{\ast}_{M/G}f$.
We have
\begin{gather*}
\nabla^A_X(fY)(x)  = T\pi_{M/G}\nabla_{X^h}(\tilde{f}Y^h)(q)= T\pi_{M/G}\big(\tilde{f}(q)\nabla_{X^h}Y^h(q)+ (\mathscr{L}_{x^h}\tilde{f})(q)\cdot Y^h(q)\big)\\
\phantom{\nabla^A_X(fY)(x)}{}= f(x)\nabla^A_XY(x)+ (\mathscr{L}_{X}f)(x)\cdot Y(x),
\end{gather*}
where the last part follows since the integral curves of $X^h$ project to integral curves of $X$. The fact that $\nabla^A$ satisf\/ies all the other properties of an af\/f\/ine connection is easily verif\/ied.
\end{proof}

\subsection{The reduced geodesic spray and its decomposition}\label{section4.2}

In this section we carry out the reduction of the geodesic spray corresponding to an invariant af\/f\/ine connections.
Let $M(M/G, G)$ be a principal f\/iber bundle with total space $M$ and structure group $G$. We denote the action of $G$ on $M$ by $\Phi: G \times M \rightarrow M$ and its tangent lift by $\Phi^T: G \times TM \rightarrow TM$. The tangent bundle projection is denoted by $\tau_M:TM \rightarrow M$.
 We can def\/ine a map $[\tau_{M}]_G: TM/G \rightarrow M/G$ as follows:
 \begin{gather*}
[\tau_M]_G([v]_G)= [\tau_M(v)]_G, \qquad [v]_G \in TM/G.
\end{gather*}
It is easy to see that $[\tau_M]_G: TM/G \rightarrow M/G$ is a vector bundle.
The adjoint bundle with $\mathfrak{g}$ as the f\/iber and $M/G$ as the base space will be represented by $\tilde{\mathfrak{g}}_{M/G}$. A typical element of $\tilde{\mathfrak{g}}_{M/G}$ will be denoted by $[x, \xi]_G$, where $x \in M$ and $\xi \in \mathfrak{g}$.
We shall also denote the tangent bundle of $M/G$ by $\tau_{M/G}: T(M/G)\rightarrow M/G$. If $A$ is a principal connection on the bundle $\pi_{M/G}: M \rightarrow M/G$, we can decompose the bundle $TM/G$ into its horizontal and vertical parts \cite{MaRa:99}.
\begin{Lemma}
\label{bundlepic}
The map $\alpha_A:TM/G \rightarrow T(M/G) \oplus \tilde{\mathfrak{g}}_{M/G}$ given by
\begin{gather*}
\alpha_A([v_x]_G)= T\pi_{M/G}(v_x)\oplus [x, A(v_x)]_G, \qquad v_x \in T_xM,\qquad x \in M
\end{gather*}
 is a vector bundle isomorphism.
\end{Lemma}
We denote the $\tilde{\mathfrak{g}}$ component of $\alpha_A$ by $\rho_A:TM/G \rightarrow \tilde{\mathfrak{g}}$.
This decomposition of $TM/G$ is $A$-dependent and we write $TM/G \simeq T(M/G)\oplus_A\tilde{\mathfrak{g}}_{M/G}$.
The lifted action $\Phi^T$ makes $TM$ the total space of a principal bundle over $TM/G$ with structure group $G$.
We denote the canonical projection by $\pi_{TM/G}:TM \rightarrow TM/G$. Furthermore, $G$ acts on $TTM$ by the tangent lift of $\Phi^T$. We denote by $\overline{T\tau_M}: TTM/G \rightarrow TM/G$ the map given by
\begin{gather*}
\overline{T\tau_M}([W_{v_x}]_G)= [T\tau_M(W_{v_x})]_G,\qquad W_{v_x} \in T_{v_x}TM.
\end{gather*}
It is easy to see that this map is well-def\/ined.
By Lemma~\ref{bundlepic}, a principal connection $\hat{A}$ on $TM(TM/G, G)$ induces an isomorphism between bundles $TTM/G$ and $T(TM/G)\oplus \tilde{\mathfrak{g}}_{TM/G}$ over $TM/G$. Thus, if $A$ and $\hat{A}$ are chosen, we can consider an identif\/ication of
$TTM/G$ and  $TT(M/G)\oplus_A T\tilde{\mathfrak{g}}_{M/G}\oplus_{\hat{A}}\tilde{\mathfrak{g}}_{TM/G}$ where we identify $T(TM/G)$ and $T(T(M/G))\oplus_{A}T\tilde{\mathfrak{g}}$ using the map $T\alpha_A$.
The following Lemma will be useful in our decomposition of the geodesic spray.
\begin{Lemma}
Given a principal connection $A$ on $M(M/G, G)$, the pullback $\hat{A}:= \tau_M^{\ast}A$ is a~principal connection on $TM(TM/G, G)$.
\end{Lemma}

 The connection $\hat{A}$ has the following useful property.
 \begin{Corollary}
If $S: TM \rightarrow TTM$ is a second-order vector field, then
\begin{gather*}
\rho_{\hat{A}}([S(v_x)]_G)= [v_x, A(v_x)]_G \in \tilde{\mathfrak{g}}_{TM/G}.
\end{gather*}
\end{Corollary}
In other words, if we choose connections $A$ and $\hat{A}$ on $M(M/G, G)$ and $TM(TM/G, G)$ respectively, studying a second-order vector f\/ield such as the geodesic spray reduces to studying the $TT(M/G)\oplus T\tilde{\mathfrak{g}}$ components, since the $\tilde{\mathfrak{g}}_{TM/G}$ component is completely determined by $A$ itself. We now def\/ine the reduced geodesic spray.\begin{Proposition}
\label{redspray}
Let $\omega$ be a $G$-invariant linear connection on $L(M)$ and $\nabla$ the corresponding connection on $M$. The map $\overline{Z}: TM/G \rightarrow TTM/G$ given by
\begin{gather*}
\overline{Z}([v_x]_G)= [Z(v_x)]_G= [T_u\Phi_{\xi}B(\xi)_u]_G
\end{gather*}
is well-defined. We call $\overline{Z}$ the  reduced geodesic spray.
\end{Proposition}
\begin{proof}
$G$-invariance of $\omega$ implies the invariance of the standard horizontal vector f\/ields. The result now follows from the $G$-equivariance of the association map $\Phi_{\xi}: L(M) \rightarrow TM$, $\xi \in \mathbb{R}^n$.
\end{proof}

Now, since $G$ acts on $TTM$ via the lifted action, we can def\/ine a map $\overline{T\pi_{TM/G}}: TTM/G \rightarrow T(T(M/G))$ as follows:
\begin{gather*}
\overline{T\pi_{TM/G}}[W_{v_x}]_G= T\pi_{TM/G}(W_{v_x}),\qquad [W_{v_x}]_G \in TTM/G.
\end{gather*}
This  is well-def\/ined since given any $g \in G$, we have
$\pi_{TM/G}\circ T\Phi^T_g= \pi_{TM/G}$. By abuse of notation, we shall use the maps $T\pi_{TM/G}$ and $\overline{T\pi_{TM/G}}$ interchangeably.
We have the following result.
\begin{Proposition}
\label{ssordvf}
Let $S_Z: T(M/G)\rightarrow TT(M/G)$ be the map defined by
\begin{gather*}
S_Z(\bar{X})= T(T\pi_{M/G}\circ \pi_{TM/G})\overline{Z}([\bar{X}^h(x)]_G),\qquad \bar{X}\in T_{[x]_G}(M/G),
\end{gather*}
where $\bar{X}^h$ is an invariant horizontal vector field that projects to $\bar{X}$ at $x \in M$.
The following statements hold:
\begin{enumerate}\itemsep=0pt
\item[$(i)$] $S_Z$ is a second-order vector field on $T(M/G)$;
\item[$(ii)$] $S_Z(\bar{X})= \bar{X}^T(\bar{X})- \mathrm{vlft}_{\bar{X}}T\pi_{M/G}(\nabla_{\bar{X}^h}\bar{X}^h)$, where, by abuse of notation, $\bar{X}$ is a vector field on $M/G$ which has a value $\bar{X}$ at $[x]_G \in M/G$.
\end{enumerate}
\end{Proposition}
\begin{proof}
$(i)$
 We compute
\begin{gather*}
T\tau_{M/G}S_Z(\bar{X}) = T(\tau_{M/G}\circ T\pi_{M/G}) Z\big(\bar{X}^h(x)\big)\\
\phantom{T\tau_{M/G}S_Z(\bar{X})}{}
 = T\pi_{M/G} T\tau_M \big(Z\big(\bar{X}^h(x)\big)\big)
 = T \pi_{M/G}\big(\bar{X}^h(x)\big)= \bar{X}.
\end{gather*}

$(ii)$ Let $\Phi^{\bar{X}^h}_t$ and $\Phi^{\bar{X}}_t$ be the f\/lows of $\bar{X}^h$ and $\bar{X}$ respectively. We have
\begin{gather*}
  TT\pi_{M/G}\big(\bar{X}^h\big)^T\big(\bar{X}^h(x)\big) = \left.\frac{d}{dt}\right|_{t= 0}\big(T\pi_{M/G}\circ T\Phi_t^{X^{h}}\big)\big(\bar{X}^h(x)\big) \\
\phantom{TT\pi_{M/G}\big(\bar{X}^h\big)^T\big(\bar{X}^h(x)\big)}{}
 =  \left.\frac{d}{dt}\right|_{t=0}T\big(\pi_{M/G}\circ \Phi_t^{X^{h}}\big)\big(\bar{X}^h(x)\big)\\
\phantom{TT\pi_{M/G}\big(\bar{X}^h\big)^T\big(\bar{X}^h(x)\big)}{}
 = \left.\frac{d}{dt}\right|_{t=0}\left.\left(\frac{d}{ds}\right|_{s=0}\big(\pi_{M/G}\circ \Phi_t^{X^{h}}\big)\big(\Phi_s^{X^{h}}(x)\big)\right)\\
\phantom{TT\pi_{M/G}\big(\bar{X}^h\big)^T\big(\bar{X}^h(x)\big)}{}
 = \left.\left.\frac{d}{dt}\right|_{t=0}\frac{d}{ds}\right|_{s= 0}\big(\pi_{M/G}\circ \Phi_{t+s}^{X^{h}}(x)\big)\\
\phantom{TT\pi_{M/G}\big(\bar{X}^h\big)^T\big(\bar{X}^h(x)\big)}{}
 = \left.\left.\frac{d}{dt}\right|_{t=0}\frac{d}{ds}\right|_{s= 0}\big(\Phi_{t+s}^{X}([x]_G\big)= \bar{X}^T(\bar{X}).
\end{gather*}

Next, we look at
\begin{gather*}
TT\pi_{M/G}\mathrm{vlft}_{\bar{X}^h(x)}\big(\nabla_{\bar{X}^h}\bar{X}^h(x)\big)
 = \left.\frac{d}{dt}\big(tT\pi_{M/G}\nabla_{\bar{X}^h}\bar{X}^h(x)+ T\pi_{M/G}\bar{X}^h(x)\big)\right|_{t= 0}\\
 \phantom{TT\pi_{M/G}\mathrm{vlft}_{\bar{X}^h(x)}\big(\nabla_{\bar{X}^h}\bar{X}^h(x)\big)}{}
   = \left.\frac{d}{dt}\big(tT\pi_{M/G}\nabla_{\bar{X}^h}\bar{X}^h(x)+ \bar{X}([x]_G)\big)\right|_{t= 0}\\
\phantom{TT\pi_{M/G}\mathrm{vlft}_{\bar{X}^h(x)}\big(\nabla_{\bar{X}^h}\bar{X}^h(x)\big)}{}
 = \mathrm{vlft}_{\bar{X}([x]_G)}T\pi_{M/G}\nabla_{\bar{X}^h}\bar{X}^h(x).
 \end{gather*}
This gives us
\begin{gather*}
S_Z(\bar{X})= \bar{X}^T(\bar{X})- \mathrm{vlft}_{\bar{X}}\big(T\pi_{M/G}\nabla_{\bar{X}^h}\bar{X}^h(x)\big).
\end{gather*}
The result now following from Proposition \ref{eq:redaff}.\end{proof}

{\sloppy The idea here is that we use principal connections $A$ and $\tau^{\ast}_MA$ on $M(M/G, G)$ and $TM(TM/G, G)$, respectively, to write the reduced geodesic spray corresponding to an invariant linear connection as a map from $T(M/G)\oplus \tilde{\mathfrak{g}}$ to $TT(M/G)\oplus T\tilde{\mathfrak{g}}$. The map $S_Z$ gives us one component of this decomposition. From Proposition \ref{zdecomp}, we see that $S_Z$ is the geodesic spray of the af\/f\/ine connection $\nabla^A$ on $M/G$.

}

Next, we def\/ine a map $P_Z: \tilde{\mathfrak{g}}\rightarrow TT(M/G)$ as follows:
\begin{gather*}
P_Z([x, \xi]_G)= TT\pi_{M/G}T\pi_{TM/G}\overline{Z}([\xi^V_L(x)]_G),
\end{gather*}
where $\xi^V_L$ is the left-invariant vector f\/ield on $M$ that satisf\/ies
$\xi^V_L(x)= \xi_M(x)$. We must verify that this is well-def\/ined. To see this, notice that $[g\cdot x, \mathrm{Ad}_g\xi]_G= [x, \xi]_G$. Next, we have
\begin{gather*}
(\mathrm{Ad}_g\xi)_M(g \cdot x) =
\left.\frac{d}{dt}\Phi_{\mathrm{exp}(\mathrm{Ad}_g\xi) t}(g \cdot x)\right|_{t= 0}
 = \left. \frac{d}{dt}\Phi\big(g(\exp \xi t)g^{-1}, g \cdot x\big)\right|_{t= 0} \\
 \phantom{(\mathrm{Ad}_g\xi)_M(g \cdot x)}{}
 = T_x\Phi_g\xi_M(x)= \xi^V_L(g \cdot x).
\end{gather*}
Let us denote $\tilde{\xi}:= [x, \xi]_G$.
Using \eqref{eqref:zxtvlf}, we get
\begin{gather*}
P_Z(\tilde{\xi})=TT\pi_{M/G}Z\big(\xi^V_L(x)\big) = TT\pi_{M/G}\big(\big(\xi^V_L\big)^T\big(\xi^V_L(x)\big)\big)- TT\pi_{M/G}\mathrm{vlft}_{\xi^V_L(x)}\big(\nabla_{\xi^V_L}\xi^V_L(x)\big)\\
\phantom{P_Z(\tilde{\xi})=TT\pi_{M/G}Z\big(\xi^V_L(x)\big)}{}
 = - \mathrm{vlft}_{{\mathbf 0}}\big(T\pi_{M/G}\nabla_{\xi^V_L}\xi^V_L(x)\big).
\end{gather*}
We write $\mathscr{S}(\tilde{\xi}, \tilde{\xi})= (T\pi_{M/G}\nabla_{\xi^V_L}\xi^V_L(x))$. Since $\nabla$ is $G$-invariant, this map is well-def\/ined.

Next, we def\/ine $R_Z: T(M/G)\rightarrow T\tilde{\mathfrak{g}}$ by
\begin{gather*}
R_Z(\bar{X})= T\rho_A T\pi_{TM/G}Z\big(\bar{X}^h(x)\big).
\end{gather*}
Then, using \eqref{eqref:zxtvlf}, we calculate
\begin{gather*}
 T\rho_A T\pi_{TM/G}Z\big(\bar{X}^h(x)\big) = T\rho_A T\pi_{TM/G}\big(\big(\bar{X}^h\big)^T\big(\bar{X}^h(x)\big) - \mathrm{vlft}_{\bar{X}^h(x)}\big(\nabla_{\bar{X}^h}\bar{X}^h(x)\big)\big).
\end{gather*}
Let us look at the f\/irst term on the right-hand side
\begin{gather*}
 T\rho_A T\pi_{TM/G}\big(\big(\bar{X}^h\big)^T\big(\bar{X}^h(x)\big)\big)  = T\rho_A T\pi_{TM/G}\left. \frac{d}{dt}\right|_{t= 0}T\Phi^{\bar{X}^h}_t\big(\bar{X}^h(x)\big)\\
 \phantom{T\rho_A T\pi_{TM/G}\big(\big(\bar{X}^h\big)^T\big(\bar{X}^h(x)\big)\big)}{}
 = \left.\frac{d}{dt}\right|_{t= 0}\rho_A \big([T\Phi^{\bar{X}^h}_t\big(\bar{X}^h(x)\big)]_G\big)= 0,
\end{gather*}
since $T\Phi^{\bar{X}^h}_t(\bar{X}^h(x))$ is horizontal and $\rho_A$ vanishes on horizontal vectors. Also,
\begin{gather*}
T\rho_AT\pi_{TM/G} \mathrm{vlft}_{\bar{X}^h(x)} \big(\nabla_{\bar{X}^h}\bar{X}^h(x)\big)
= \left.\frac{d}{dt}\right|_{t= 0}\rho_A \pi_{TM/G}\big(t\nabla_{\bar{X}^h}\bar{X}^h(x)+ \bar{X}^h(x)\big)\\
\phantom{T\rho_AT\pi_{TM/G} \mathrm{vlft}_{\bar{X}^h(x)}}{}
 = \left.\frac{d}{dt}\right|_{t= 0}\big(t\rho_A\circ \pi_{TM/G}\big(\nabla_{\bar{X}^h}\bar{X}^h(x)\big)+ \rho_A\circ \pi_{TM/G}\big(\bar{X}^h(x)\big)\big)\\
\phantom{T\rho_AT\pi_{TM/G} \mathrm{vlft}_{\bar{X}^h(x)}}{}
 = \left.\frac{d}{dt}\right|_{t= 0}\big(t\rho_A\circ \pi_{TM/G}\big(\nabla_{\bar{X}^h}\bar{X}^h(x)\big)+ 0\big),
\end{gather*}
and thus we get
\begin{gather*}
R_Z(\bar{X})
= -\mathrm{vlft}_{0}\big(\rho_A\circ \pi_{TM/G}\big(\nabla_{\bar{X}^h}\bar{X}^h(x)\big)\big).
\end{gather*}
If $HM$ is geodesically invariant, then $\nabla_{\bar{X}^h}\bar{X}^h$ is horizontal, and thus $R_Z= 0$.

Finally, we def\/ine $U_Z: \tilde{\mathfrak{g}}\rightarrow T\tilde{\mathfrak{g}}$ by
\begin{gather*}
U_Z(\tilde{\xi})= T\rho_AT\pi_{TM/G}Z\big(\xi^V_L(x)\big),
\end{gather*}
and a calculation similar to the one performed above shows that
\begin{gather*}
U_Z(\tilde{\xi})= -\mathrm{vlft}_{\xi}\rho_A\big(\pi_{TM/G}\big(\nabla_{\xi^V_L}\xi^V_L(x)\big)\big).
\end{gather*}
The following lemma is useful.
\begin{Lemma}
The map  $\widetilde{\nabla}^A: \Gamma(T(M/G)) \times \Gamma(\tilde{\mathfrak{g}})\rightarrow \Gamma(\tilde{\mathfrak{g}})$ given by
\begin{gather*}
\widetilde{\nabla}^A_{\bar{X}}\tilde{\xi}([x]_G)= \rho_A\pi_{TM/G}\big(\big\langle \bar{X}^h: \xi^V_L\big\rangle (x)\big), \qquad [x]_G \in (M/G)
\end{gather*}
defines a vector bundle connection on the bundle $\tilde{\mathfrak{g}}$.
\end{Lemma}
\begin{proof}
Let $f: M/G \rightarrow \mathbb{R}$ be a dif\/ferentiable function. Def\/ine $f^h: M \rightarrow \mathbb{R}$ by $f^h= \pi_{M/G}^{\ast}f$. Therefore, $(f\bar{X})^h= f^h\bar{X}^h$. We compute
\begin{gather*}
\widetilde{\nabla}^A_{f\bar{X}}\tilde{\xi}  = \rho_A\pi_{TM/G}\big(\big\langle f^h\bar{X}^h: \xi^V_L\big\rangle\big)\\
\phantom{\widetilde{\nabla}^A_{f\bar{X}}\tilde{\xi}}{}
 = \rho_A\pi_{TM/G}\big(f^h\nabla_{\bar{X}^h}\xi^V_L+ f^h\nabla_{\xi^V_L}\bar{X}^h+ \big(\mathscr{L}_{\xi^V_L}f^h\big)\bar{X}^h \big)
 = f\widetilde{\nabla}^A_{\bar{X}}\tilde{\xi},
\end{gather*}
since $(\mathscr{L}_{\xi^V_L}f^h)\bar{X}^h= 0$. The property $\widetilde{\nabla}^A_{\bar{X}}f\tilde{\xi}= f\widetilde{\nabla}^A_{\bar{X}}\tilde{\xi}+ (\mathscr{L}_{\bar{X}}f)\tilde{\xi}$
can be proved similarly.
\end{proof}

We now state the main result of this section.
\begin{Theorem}
\label{main}
Let $Z_h: T(M/G)\oplus \tilde{\mathfrak{g}} \rightarrow TT(M/G)$ be the map defined by
\begin{gather*}
Z_h (\bar{X}\oplus \tilde{\xi} )= TT\pi_{M/G}\overline{Z}\big[\bar{X}^h(x)+ \xi^V_L(x)\big]_G,
\end{gather*}
where  $\bar{X}^h$ is an invariant horizontal vector field that projects to $\bar{X}$ at $x \in M$, and $\xi^V_L$ is the left-invariant vertical vector field with value $\xi_M(x)$ at $x \in M$.

Let $Z_v:T(M/G)\oplus \tilde{\mathfrak{g}} \rightarrow T\tilde{\mathfrak{g}}$ be the map defined by
\begin{gather*}
Z_v(\bar{X}\oplus \tilde{\xi})= T\rho_A\overline{Z}\big(\big[\bar{X}^h(x)+ \xi^V_L(x)\big]_G\big),
\end{gather*}
where $\bar{X}^h$ and $\xi^V_L$ are defined as above.
The following statements hold:
\begin{enumerate}
\item[$(i)$] $
Z_h(\bar{X}\oplus \tilde{\xi})= S_Z(\bar{X})- \mathrm{vlft}_{\bar{X}}\mathscr{S}(\tilde{\xi}, \tilde{\xi})- \mathrm{vlft}_{\bar{X}}\big(T\pi_{M/G}\big\langle \bar{X}^h: \xi^V_L \big\rangle \big) $;
\item[$(ii)$]
$Z_v(\tilde{X}\oplus \tilde{\xi})
= R_Z(\bar{X})+ U_Z(\tilde{\xi})- \mathrm{vlft}_{\xi}\big(\widetilde{\nabla}^A_{\bar{X}}\tilde{\xi}([x]_G) \big)$.
\end{enumerate}
\end{Theorem}

\begin{proof}
Let us compute
\begin{gather*}
TT\pi_{M/G} Z\big(\bar{X}^h(x) + \xi^V_L(x)\big)\\
\phantom{TT\pi_{M/G}}{}
=  TT\pi_{M/G}\big(\big(\bar{X}_h+ \xi^V_L\big)^T\big(\bar{X}^h+ \xi^V_L(x)\big)\big)- \mathrm{vlft}_{\bar{X}^h(x)}\big(T\pi_{M/G}\nabla_{\bar{X}^h+ \xi^V_L}(X^h+ \xi^V_L)\big)\\
\phantom{TT\pi_{M/G}}{}
=  TT\pi_{M/G}\big(\bar{X}^h\big)^T(\xi_M(x))+ \bar{X}^T(\bar{X}([x]_G)- \mathrm{vlft}_{\bar{X}^h(x)}\big(T\pi_{M/G}\nabla_{\bar{X}^h}\bar{X}^h\big)\\
\phantom{TT\pi_{M/G}=}{}
-\mathrm{vlft}_{\bar{X}^h(x)}\big(\mathscr{S}(\tilde{\xi}, \tilde{\xi})\big)- \mathrm{vlft}_{\bar{X}^h(x)}\big(T\pi_{M/G}\big\langle \bar{X}^h: \xi^V_L \big\rangle\big)\\
\phantom{TT\pi_{M/G}}{}
=   TT\pi_{M/G}\big(\bar{X}^h\big)^T\big(\xi^V_L(x)\big)+ S_Z(\bar{X})- \mathrm{vlft}_{\bar{X}^h(x)}\big(\mathscr{S}(\tilde{\xi}, \tilde{\xi})\big)\\
\phantom{TT\pi_{M/G}=}{}
 - \mathrm{vlft}_{\bar{X}^h(x)}\big(T\pi_{M/G}\big\langle \bar{X}^h: \xi^V_L\big\rangle(x)\big).
\end{gather*}
We also have
\begin{gather*}
TT\pi_{M/G}\big(\bar{X}^h+ \xi^V_L\big)^T\big(\xi^V_L(x)\big)= \left.\frac{d}{dt}\right|_{t= 0}T\pi_{M/G} \xi^V_L\big(\Phi^{\bar{X}^h}_t(x)\big)= 0.
\end{gather*}
This gives us the f\/irst part. Part $(ii)$ follows from a similar computation.
\end{proof}

\begin{Remark}
The fact that the right-hand sides of $Z_h$ and $Z_v$ respectively are independent of the extensions follows from $G$-invariance of $\omega$ and the def\/inition of $\overline{Z}$.
\end{Remark}
\begin{Remark}
The decomposition of the reduced geodesic spray into horizontal and vertical parts given in Theorem~\ref{main} is similar to the decomposition of second-order systems in Crampin and Mestdag~\cite{crme:09}, particularly in the case of an af\/f\/ine spray. We discuss the Riemannian case in Section~\ref{section5.1}.
\end{Remark}

\section{Discussion}\label{section5}

The horizontal part of the reduced geodesic spray consists of three terms. The map $S_Z$ is a~second-order vector f\/ield on $T(M/G)$. The term $\mathscr{S}(\tilde{\xi}, \tilde{\xi})$ can be interpreted in the following manner. Recall that the second fundamental form corresponding to the vertical distribution is a map $\mathcal{S}: \Gamma(VM) \times \Gamma (VM) \rightarrow HM$ def\/ined by
\begin{gather*}
\mathcal{S}(v_x, w_x)= \mathrm{hor}\left(\nabla_XY\right),\qquad v_x,  w_x \in V_xM,
\end{gather*}
where $X$ and $Y$ are extensions of $v_x$ and $w_x$ respectively.
In view of this, we have
\begin{gather*}
\mathcal{S}(\xi_M(x), \xi_M(x))= \big(\mathscr{S}(\tilde{\xi}, \tilde{\xi})\big)^h(x).
\end{gather*}
Now, the vertical distribution $VM$ is geodesically invariant if and only if $\mathcal{S}$ is skew-symmetric. Hence, if $VM$ is geodesically invariant, we have $\mathscr{S}(\tilde{\xi}, \tilde{\xi})= 0$.

\subsection{The Riemannian case}\label{section5.1}

The last term in the horizontal part of the reduced geodesic spray is related to the curvature of the horizontal distribution, at least in the case when $M$ is a Riemannian manifold with an invariant Riemannian metric, the chosen af\/f\/ine connection is the Levi-Civita connection corresponding to this metric, and $A$ is the mechanical connection as we show below.

Let $(M, k)$ be a Riemannian manifold and $G$ be a Lie group that acts freely and properly on~$G$, so that $\pi_{M/G}: M \rightarrow M/G$ is a principal bundle. Suppose that the Riemannian metric $k$ is invariant under $G$. The mechanical connection corresponding to $k$ is a principal connection on $\pi_{M/G}:M \rightarrow M/G$ determined by the condition that the horizontal subbundle is orthogonal to the vertical subbundle $VM$ with respect to the metric. We denote by $A$ the connection one-form corresponding to this connection. We also let $\nabla$ be the Levi-Civita connection corresponding to~$k$.

\begin{Lemma}
The following holds
\begin{gather*}
k\big(\big\langle \bar{X}^h: \xi^V_L\big\rangle (x), \bar{Y}^h(x)\big)= k\big(\big(B_A\big(\bar{X}^h(x), \bar{Y}^h(x)\big)\big)_M, \xi^V_L(x)\big),
\end{gather*}
where $\bar{X}^h$ and $\bar{Y}^h$ are invariant horizontal vector fields on $M$, and $B_A$ is the curvature form corresponding to $A$.
\end{Lemma}

\begin{proof}
Recall that if $X$, $Y$ and $Z$ are vector f\/ields on $M$, the Koszul formula is given by
\begin{gather*}
  2k(\nabla_X Y,Z)  = \mathscr{L}_X(k(Y, Z))+ \mathscr{L}_Y(k(X, Z))- \mathscr{L}_Z(k(X, Y))+ k([X, Y], Z)\\
  \phantom{2k(\nabla_X Y,Z)  =}{}
  - k([X, Z], Y])- k([Y, Z], X).
\end{gather*}
We therefore have (using the Koszul formula twice and adding the two results)
\begin{gather*}
 2k\big(\big\langle \bar{X}^h:\xi^V_L\big\rangle(x), \bar{Y}^h(x)\big)
 = 2\mathscr{L}_{\bar{X}^h}\big(k\big(\bar{Y}^h(x), \xi^V_L(x)\big)\big)
 + 2\mathscr{L}_{\xi^V_L}\big(k\big(\bar{X}^h(x), \bar{Y}^h(x)\big)\big)\\
\qquad{} - 2\mathscr{L}_{\bar{Y}^h}\big(k\big(\bar{X}^h(x), \xi^V_L(x)\big)\big)
- 2k\big(\big[\bar{X}^h, \bar{Y}^h\big](x), \xi^V_L(x)\big)  - 2k\big(\big[\xi^V_L, \bar{Y}^h\big](x), \bar{X}^h(x)\big).
\end{gather*}
Now, the f\/irst and the third terms respectively on the right-hand side are clearly zero (by the def\/inition of the mechanical connection). The second term is zero since the function $k\left(\bar{x}^h(x), \bar{Y}^h(x)\right)$ is constant along the invariant vertical vector f\/ield $\xi^V_L$. The f\/ifth term is also zero since the Lie bracket $[\xi^V_L, \bar{Y}^h]$ is a vertical vector f\/ield. Thus, we get
\begin{gather*}
k\big(\big\langle \bar{X}^h: \xi^V_L\big\rangle (x), \bar{Y}^h(x)\big)= k\big(\big[\bar{X}^h, \bar{Y}^h\big](x), \xi^V_L(x)\big).
\end{gather*}
By the Cartan structure formula, we have
\begin{gather*}
\big[\bar{X}^h, \bar{Y}^h\big]=  [\bar{X}, \bar{Y}]^h- \big(B_A\big(\bar{X}^h, \bar{Y}^h\big)\big)_M(x).
\end{gather*}
Therefore,
\begin{gather*}
k\big(\big\langle \bar{X}^h: \xi^V_L\big\rangle (x), \bar{Y}^h(x)\big)= k\big(\big(B_A\big(\bar{X}^h(x), \bar{Y}^h(x)\big)\big)_M, \xi^V_L(x)\big).\tag*{\qed}
\end{gather*}
 \renewcommand{\qed}{}
\end{proof}

The vertical part of the reduced geodesic spray consists of the map $R_Z$ which vanishes identically if the horizontal distribution corresponding to the principal connection $A$ is geodesically invariant, and can be thought of as the fundamental form corresponding to the horizontal distribution. Lewis \cite{ADL:98} has shown that if both $HM$ and $VM$ are geodesically invariant, then the corresponding linear connection restricts to the subbundle $L(M, A)$.
The term $U_Z(\tilde{\xi})$ is essentially the Euler\textendash{}Poincar\'e term, and the last term corresponds to a connection on $\tilde{\mathfrak{g}}$.

\subsection{Forces in mechanics}\label{section5.2}
As mentioned in the Introduction, our motivation for studying reduction in the af\/f\/ine connection setup comes from mechanics. In this sense, studying the geodesic spray corresponds to looking at mechanical systems with no external forces (in other words, the dynamics are given by the geodesic equation). It is worth considering the case in which forces are present as many important examples in mechanics fall in this class. In the following, we consider the so-called simple mechanical systems \cite{AbMa:78, BuLe:05, Sma:70}. A simple mechanical system is a triple $(M,k, V)$ where $(M, k)$ is Riemannian manifold and $V: M \rightarrow \mathbb{R}$ is a smooth function (called the potential function). The gradient of $V$ is a vector f\/ield on $M$ def\/ined by
\begin{gather*}
\mathrm{grad}V(x)= k^{\#}(dV(x)),
\end{gather*}
where $k^{\#}: T^{\ast}M\rightarrow TM$ is a vector bundle isomorphism over $M$ induced by the metric $k$.
The dynamics of such a system are given by
\begin{gather}
\label{eq:forces}
\nabla_{c^{\prime}(t)}c^{\prime}(t)=- (\mathrm{grad} V)(c(t)),\qquad c(t) \in M,
\end{gather}
where $\nabla$ is the Levi-Civita connection on $M$.
Equivalently, one can study the following equation on $TM$ \cite{AbMa:78}:
\begin{gather}
\label{eq:tmforces}
v^{\prime}(t)= Z(v(t))- \mathrm{vlft}_{v(t)}\left(\mathrm{grad} V(c(t))\right),\qquad v(t) \in TM,\qquad \tau_M(v(t))= c(t).
\end{gather}
In the unforced case dealt with in this paper, we study the geodesic spray because its integral curves project to geodesics on $M$. Similarly, we can study the second-order vector f\/ield
\begin{gather*}
\widetilde{Z}(v)= Z(v)- \mathrm{vlft}_{v}\left(\mathrm{grad} V(\tau_M(v))\right),
\end{gather*}
which, by Proposition~\ref{zdecomp}, can be written as
\begin{gather*}
\widetilde{Z}(v)= X^{T}(v)- \mathrm{vlft}_v\left(\nabla_{X_v}X_v+ \mathrm{grad} V(\tau_M(v))\right).
\end{gather*}
If $\nabla$, $k$ and $V$ are $G$-invariant, we can study the reduction of $\widetilde{Z}$ using our methodology. Even though we consider a potential force here, a general force $F$ can be incorporated in this picture by essentially replacing $\mathrm{grad} V$ with $k^{\#}(F)$ in~\eqref{eq:forces}.

\subsection{Generalized connections}\label{section5.3}

In our investigation, we have considered a $G$-invariant af\/f\/ine connection $\nabla$ on a manifold $M$ along with a principal connection $A$ on $\pi_{M/G}:M \rightarrow M/G$ and used it to def\/ine an $A$-dependent af\/f\/ine connection $\nabla^A$ on $M/G$. Equivalently, and perhaps more naturally, $\nabla$ induces a connection on the vector bundle $[\tau_M]_G:TM/G\rightarrow M/G$. Since the principal connection $A$ provides a~decomposition of $TM/G$, we can also recover $\nabla^A$ in this manner.

Furthermore, we can consider a generalized connection \cite{CaLa:03,CoMa:04} on the vector bundle $[\tau_M]_G: TM/G \rightarrow M/G$ and explore how our reduction procedure can be applied to this more general situation.

\subsection{Nonholonomic systems with symmetry}\label{section5.4}
Roughly speaking, nonholonomic systems are mechanical systems with velocities constrained to lie in a given non-integrable distribution.
Following the fundamental paper of Koiller \cite{Koi:92} there has been a lot of interest in
studying symmetries and reduction of nonholonomic systems \cite{Bat:02,BKMM:96, BuLe:05, CaLeMaDi:98, Cor:02, ADL:98, SnJp:98}. In this section we outline how these systems can be studied in our framework.
Let $(M, k)$ be a Riemannian manifold with a $G$-invariant Riemannian metric. Let $\mathcal{D}$ be a smooth, non-integrable, $G$-invariant distribution on $M$ and $\mathcal{D}^{\perp}$ the orthogonal complement with respect to the metric $k$. Let $\nabla$ be the Levi-Civita af\/f\/ine connection associated with $k$. The Lagrange--d'Alembert
principle allows us to conclude that the constrained geodesics $c(t) \in
M$ satisfy \cite{ADL:98}
\begin{gather*}
{\nabla}_{c^{\prime}(t)}c^{\prime}(t) \in
\mathcal{D}_{c(t)}^{\perp},\qquad   c^{\prime}(t) \in
\mathcal{D}_{c(t)}.
\end{gather*}
Sometimes these conditions are written as
\begin{gather*}
{\nabla}_{c^{\prime}(t)}c^{\prime}(t) = \lambda(c(t)), \qquad
P^{\perp}(c^{\prime}(t))  =0,
\end{gather*}
where $\lambda$ is a section of $\mathcal{D}^{\perp}$ and $P^{\perp}: TM
\rightarrow TM$ is
the projection onto $\mathcal{D}^{\perp}$. It can be shown that the trajectories $c: \mathbb{R} \rightarrow M$ satisfying the constraints are actually geodesics of an af\/f\/ine
connection
$\widetilde{\nabla}$ def\/ined by $\widetilde{\nabla}_XY= {\nabla}_XY
+ (\nabla_XP^{\perp})(Y)$.
Note that, in general, the connection $\widetilde{\nabla}$ (sometimes called a constrained connection) will not be Levi-Civita.
We can use our approach to study the geodesic spray of the constrained connection in the presence of a principal connection on
 $\pi_{M/G}: M \rightarrow M/G$.
In such a case, the picture gets more complicated since a decomposition of $TTM/G$, and that of the geodesic spray, will depend on the distribution $\mathcal{D}$, and we hope to address this problem in subsequent work.
This procedure is related to the reduction of ``external'' symmetries of a  generalized $G$-Chaplygin system~\cite{EKMR:05}.

\subsubsection*{Acknowledgements}

I would like to thank my thesis supervisor Dr. Andrew Lewis for his constant guidance and support. This work would not have materialized without the  many invaluable discussions I~have had with him over the years. The author also thanks the anonymous referees for their constructive comments on a previous version of this paper.

\pdfbookmark[1]{References}{ref}
\LastPageEnding

\end{document}